\documentclass[12pt]{amsart}
\newfont{\sheaf}{eusm10 scaled\magstep1}

\usepackage[all]{xy}

\newcommand{\ra}{\ensuremath{\rightarrow}}

\def\eea{\end{eqnarray*}}
\def\bea{\begin{eqnarray*}}
\def\Bbb{\bf}
\def\P{{\Bbb P}}

\def\X{{\mathcal{X}}}
\def\K{{\mathcal{K}}}
\def\de{{\delta}}
\def\De{{\Delta}}
\def\ga{{\gamma}}
\def\Ga{{\Gamma}}

\newcommand{\Proof}{{\it Proof. }}
\newcommand{\QED}{{\hfill $Q.E.D.$}}

\newtheorem{teo}{Theorem}[section]
\newtheorem{df}[teo]{Definition}
\newtheorem{lem}[teo]{Lemma}
\newtheorem{cor}[teo]{Corollary}

\newtheorem{oss}[teo]{Remark}
\newtheorem{prop}[teo]{Proposition}

\newcommand{\C}{\ensuremath{\mathbb{C}}}
\newcommand{\R}{\ensuremath{\mathbb{R}}}
\newcommand{\Z}{\ensuremath{\mathbb{Z}}}
\newcommand{\Q}{\ensuremath{\mathbb{Q}}}

\newcommand{\N}{\ensuremath{\mathbb{N}}}
\newcommand{\hol}{\ensuremath{\mathcal{O}}}
\newcommand{\HH}{\ensuremath{\mathbb{H}}}
\newcommand{\PP}{\ensuremath{\mathbb{P}}}
\newcommand{\RR}{\ensuremath{\mathcal{R}}}
\newcommand{\BB}{\ensuremath{\mathcal{B}}}
\newcommand{\FF}{\ensuremath{\mathcal{F}}}
\newcommand{\A}{\ensuremath{\mathcal{A}}}

\newcommand{\HHH}{\ensuremath{\mathcal{H}}}

\newcommand{\I}{\ensuremath{\mathcal{I}}}
\newcommand{\SSS}{\ensuremath{\mathcal{S}}}
\newcommand{\PPP}{\ensuremath{\mathcal{P}}}

\begin{document}

\title{Q.E.D. for Algebraic Varieties.}

\author{Fabrizio Catanese\\
    Universit\"at Bayreuth}

\date{First Version: September 11, 2002-Final version:  September 27, 2005.}
\maketitle

\begin{abstract}
    We introduce a new equivalence relation for complete algebraic
varieties with canonical
singularities,  generated by birational
    equivalence, by flat algebraic deformations
(of varieties with canonical singularities), and by quasi-\'etale morphisms,
i.e., morphisms which are unramified in codimension $1$. We denote
the above equivalence by A.Q.E.D. : = Algebraic-Quasi-\'Etale- Deformation.

    A completely similar equivalence relation, denoted by
  $\C$-Q.E.D.  (: = Complex-Quasi-\'Etale- Deformation),
  can be considered for
compact complex manifolds and for compact complex spaces with
canonical singularities.
It is generated by  bimeromorphic
    equivalence, by flat deformations
of  complex spaces with canonical singularities,
    and by quasi-\'etale morphisms.

By a recent theorem of Siu (conjecturally holding also for compact
K\"ahler manifolds),
dimension and Kodaira dimension are invariants for
$A.Q.E.D.$ (also for $\C$-Q.E.D. -equivalence of algebraic varieties?).

It was an interesting question whether conversely
two algebraic varieties of the same dimension and with the same
    Kodaira dimension are Q.E.D. -
equivalent (the question then bifurcates into: A.Q.E.D., or at least
$\C$-Q.E.D.).

The answer to the A.Q.E.D. question is positive for curves by well
known results.
    Using Enriques'  classification we show first that the answer  to
the $\C$-Q.E.D. question
is positive  for special algebraic surfaces
(i.e., for algebraic surfaces with Kodaira dimension at most $1$).
More generally, using Kodaira's classification we show that the same
result holds
for compact complex surfaces with
Kodaira dimension  $0, 1$ and even first Betti number.

The appendix by S\"onke Rollenske shows that the same does not hold
if we consider also surfaces with odd first Betti number. Indeed he
proves that any surface which is  $\C$-Q.E.D.-equivalent to a
Kodaira surface is itself a
Kodaira surface.

Concerning the A.Q.E.D. question, we show that the answer is positive
for complex algebraic
surfaces of Kodaira dimension $\leq 1$.

The answer to the Q.E.D. question is instead
negative for surfaces of general type:
the other appendix, due to Fritz Grunewald, is devoted to showing that
the (rigid) Kuga-Shavel type surfaces of general type obtained  as quotients of the
bidisk via discrete groups constructed from quaternion algebras
belong to countably many  distinct Q.E.D. equivalence classes.

  We leave aside here the
A.Q.E.D. question
for algebraic surfaces of Kodaira dimension $ \leq 1$ defined over an
algebraically closed
field of positive characteristic;
we hope to be able to address this question in the  future.

We  pose then several questions, and point out possible applications
of the above ideas, with the aim of showing the importance
for classification theory and for other purposes
of  investigating   Q.E.D. equivalence.

     \end{abstract}

\vfill
\pagebreak
\section{Introduction}

The purpose of the present article is to define some new and broad
equivalence relations,
called Q.E.D. ,
in the classification theory of algebraic varieties and compact complex spaces,
and to pose some  problems concerning   invariants for Q.E.D. equivalences.

To  briefly explain the prehistory of the question,
let me first recall that, in order to make the study of algebraic
varieties possible, it is customary to
introduce some equivalence relation.
The most classical one is the so called birational  equivalence,
which allows   in
particular not to distinguish between the
different projective embeddings of the same
variety (respectively, one considers the bimeromorphic equivalence
for complex spaces)

Since moreover quite often the construction of algebraic varieties
    depends upon parameters, Kodaira and Spencer introduced the notion of
{\bf $\C$-deformation equivalence} for complex manifolds:
they
\cite {k-s58}) defined two complex manifolds $X'$, $X$ to be
{\bf directly deformation equivalent} if there is a  proper holomorphic
    submersion
$ \pi : \Xi \ra
\Delta$ of a complex manifold $\Xi$ to the unit disk in the complex plane,
    such that
$X, X'$ occur as fibres of $\pi$.
If we take the equivalence relation generated
by direct deformation equivalence, we obtain the relation of
complex deformation
equivalence, and we say that $X$ is a complex
{\bf deformation} of $X'$ {\bf in the large} if
$X, X'$ are complex deformation equivalent.

These two notions extend to the case of compact complex manifolds
the classical
notions of irreducible, resp. connected, components of moduli spaces.

My first main motivation for introducing  Q.E.D.-equivalence is the
following: to explain the Kodaira classification in the case of algebraic
curves, one has  to  say that a curve has Kodaira dimension $1$ iff
it has genus $ g \geq 2$, and then to recall that  all curves of
a fixed genus $g$ are deformation equivalent.

The simple but key observation is
that for each $ g
\geq 2$ there is a curve of genus $g$ which is an \'etale
(unramified) covering
of a curve of genus $2$. Therefore all the curves with Kodaira
dimension $1$ are equivalent by the equivalence relation generated
by deformation and by \'etale maps (the same holds of course for curves
of  Kodaira dimension $0$, resp. $- \infty$).

More remarkable is the consideration of
    the (complex) algebraic surfaces
of Kodaira dimension $0$: the Enriques surfaces have an \'etale double cover
which is a $K3$-surface, and the Hyperelliptic surfaces have an \'etale cover
which is a torus (indeed, the product of two elliptic curves). So, in
this case, we should link tori and $K3$'s
by \'etale maps and deformations. This is obviously not possible, because
tori are $K(\pi,1)$'s while $K3$'s are simply connected.

That's why the solution is to divide the torus by multiplication by
$-1$, obtaining the (singular) Kummer surfaces, and then take a smoothing
of the Kummer surface to obtain a smooth $K3$ surface.
The small price to pay
is to allow morphisms
which are not \'etale, but only \'etale in codimension $1$, and moreover
to allow singularities, ordinary double points in the case of surfaces.
    This singularity is
a very special case of the Rational Double Points,
    which are the canonical singularities of dimension $2$ (cf.
\cite{reid1}. \cite{reid2}).

This remark helps to justify the following definition

\begin{df}
The relation of Algebraic Quasi-\'etale deformation is the
equivalence relation,
for complete Algebraic Varieties with canonical singularities defined
over a fixed
algebraically closed field,
generated by

\begin{itemize}
\item
(1) Birational equivalence
\item
(2) Flat proper algebraic deformation $\pi : \X \ra B$ with  Base $B$
a connected
algebraic variety, and with all the fibres
having canonical singularities
\item
(3) Quasi-\'etale morphisms $f : X \ra Y$, i.e., surjective morphisms
which are \'etale
in codimension $1$ on $X$ (there is $Z \subset X$ of codimension
$\geq 2$
such that $f$ is \'etale on $ X - Z$)

\end{itemize}

and denoted by $A.Q.E.D.$  ( $X \sim_{A.Q.E.D.} Y $).

Note that we have a completely analogous $\C-   Q.E.D.$ -equivalence
for compact complex spaces with canonical singularities generated by

\begin{itemize}
\item
(1) Bimeromorphic equivalence
\item
(2) Flat proper complex deformations $\pi : \X \ra B$ with  connected Base $B$,
    and with all the fibres
having canonical singularities
\item
(3) Quasi-\'etale morphisms $f : X \ra Y$
\end{itemize}

One may of course restrict the latter equivalence relation to algebraic
varieties, and to K\"ahler manifolds and spaces.

Finally, we define a compact complex space to be {\bf standard} if it is
$\C-   Q.E.D.$ -equivalent to a product of curves, and similarly define
the concept of an {\bf A.standard = algebraically standard} algebraic variety.
\end{df}

\begin{oss}
By Siu's recent Theorem (\cite{siu}), not only dimension, but also
the Kodaira dimension is an invariant for $A.Q.E.D.$ -equivalence if
we restrict
ourselves to consider projective varieties
with canonical singularities (defined over $\C$). It is conjectured
(ibidem, cf. also \cite{siu2}) that the
deformation invariance of plurigenera should be true more generally
for  K\"ahler complex spaces (with canonical singularities).

\end{oss}

{\bf Question 0:} Is Kodaira dimension also an invariant for $\C$-Q.E.D.
equivalence of algebraic varieties and compact K\"ahler manifolds?

We  begin with the following two theorems

\begin{teo} \label{kaehler}
Let $S$ and $S'$ be smooth K\"ahler surfaces which have the same
Kodaira dimension $\K \leq 1$. Then $S$ and $S'$ are $ \C -
Q.E.D.$-equivalent.
\end{teo}

\begin{cor} \label{0,1}
Let $S$ and $S'$ be smooth compact complex surfaces
with even first Betti numbers  and which
have the same Kodaira dimension $\K = 0, 1$. Then $S$ and $S'$ are $ \C -
Q.E.D.$-equivalent.
\end{cor}

The next theorem is not a special case of the previous, since we
only consider projective deformations

\begin{teo}\label{algebraic}
Let $S$ and $S'$ be smooth complex algebraic surfaces which have the same
Kodaira dimension $\K \leq 1$. Then $S$ and $S'$ are A.Q.E.D.-equivalent.
\end{teo}

The ingredients of the  proof of \ref{kaehler} are, beyond the
Enriques classification
and a detailed knowledge of the deformation types of elliptic surfaces,
the notion of orbifold fundamental group of a fibration, and the
following very
simple devise.

{\bf Main Observation:} Assume that we have  two effective actions
of a finite group $G$ on algebraic varieties $X$, resp. $Y$
(effective means that no element $ g \in G$ acts trivially).
Then the product action of $G$ on the product $ X \times Y$
yields a quasi - \'etale map $ X \times Y \ra (X \times Y) / G.$

A couple of words concerning A.Q.E.D.
equivalence:  in the case of Kodaira dimension $1$ we have to face
the problem that
    algebraic deformation is not completely understood
for elliptic surfaces, and, even more, the determination of
quasi-\'etale maps on
models with canonical singularities requires a rather deep understanding
of the configuration of curves allowed on some deformation of a given elliptic
surface. Such a study, as is the case for K3 surfaces, is related to a
systematic investigation of the period map for elliptic surfaces
(this investigation
was  started for Jacobian elliptic surfaces,
    i.e., elliptic surfaces admitting a section, by Chakiris in
\cite{chak1,chak2}).

Our solution to prove the Q.E.D. statement for Kodaira
dimension $1$ is to try to reduce to the case of no multiple fibres:
this is done via
the orbifold fundamental group of a fibration and works easily except
for elliptic surfaces over $\PP^1$
with one or two multiple fibres. For these,  the simplest approach
(deformation to constant moduli)
fails to work.
In this case however, the result  follows by showing the existence
    of an algebraic  deformation of such a surface to another one
possessing two (resp. : one)
$\tilde{D}_4$ fibres:
after contracting the non central $-2$ curves we get a singular
surface with a large
orbifold fundamental group, and we again reduce to the case of no
multiple fibres.

A main purpose of this article is also to pose the following

{\bf Main Question:}  which are the Q.E.D. equivalence classes
  of surfaces of general
type and of special varieties in higher dimension ?

Let us try to separate the two issues. Our proof that Kodaira
dimension and $\C$-Q.E.D. equivalence coincide for special surfaces
is somehow
related to the Def = Diff problem ( cf. \cite{f-m1}, but compare also
\cite{f-m}, 205-208). We know that two
special surfaces  $S_1, S_2$
are orientedly diffeomorphic if and only if either they
are deformation equivalent, or $S_1$ is deformation equivalent
to the complex conjugate of $S_2$.

It was recently shown ( \cite{man01}, \cite{k-k01},\cite{cat01},
\cite{c-w})
that  for algebraic surfaces of general type (for these, $\C$-Q.E.D.
and A.Q.E.D.
equivalence coincide) diffeomorphism or symplectic equivalence are no
sufficient criteria  to guarantee
    complex deformation equivalence. Moreover, we observe that almost all
(cf. Question 8)
the known counterexamples are known to be in the "standard" Q.E.D.
equivalence class, i.e.,
distinct connected components of the moduli
space of surfaces of general type are simply obtained via surfaces
which are Q.E.D. -equivalent to products of curves.

This observation leads to the second main motivation for 
introducing Q.E.D. equivalence: 
quasi- \'etale maps have for long time been an ace of diamonds 
in the sleeves of algebraic geometers in order to produce
very interesting counterexamples
(we shall point out other examples  later). Our philosophy here is that
quasi- \'etale maps are a fact of life which in classification
theory should be considered more as the daily rule rather 
than the exception.

The main questions that we want to pose can then be summarized as:

{\bf Question 1:} are there more (effectively computable) invariants for
    Q.E.D. equivalence, than  dimension and  Kodaira dimension ?

{\bf Question 2:} Is it possible to determine the  Q.E.D. equivalence
classes inside the
class of varieties with fixed dimension $n$,
and with Kodaira dimension $k$?

For curves and special algebraic surfaces over $\C$ we saw that there
is only one  A.Q.E.D. class,
but in the first appendix   Fritz Grunewald shows, considering
some Kuga-Shavel type surfaces of general type constructed
from quaternion algebras  according to general 
lines suggested
by Shimura (cf. \cite{shav})

{\bf Corollary \ref{aleph}}
{\em There are infinitely many Q.E.D. -equivalence classes
of algebraic surfaces of general type.}

The above surfaces are rigid, but the Q.E.D. -equivalence class
contains countably many distinct birational classes. We can then pose
a more daring

{\bf Question 3:} are there for instance varieties which are isolated in their
Q.E.D.-equivalence class (up to birational equivalence, of course)?

\begin{oss}

\begin{itemize}
\item
Singularities play here an essential role. Note first of all that
(as we will show in the next section), without the
restriction on these given in (2), we obtain the trivial
    equivalence relation
for algebraic varieties of the same dimension
(does this also  hold for compact K\"ahler manifolds?).
\item
Assume that a variety $X$ has the following properties of being

1) rigid,

2) smooth with ample canonical bundle,

3)  with a trivial algebraic fundamental group

4) with a trivial  group of Automorphisms.

Then any variety $X'$ birational to $X$ and with canonical singularities
has $X$ as canonical model, and since $X$ has no deformations, and there is no non
trivial quasi-\'etale map
$Y \ra X$.

The only possibility, to avoid that $X$ be isolated in its
QED-equivalence class, would be that there exists a
quasi-\'etale map $f : X \ra Y$.

If $f$ is not birational, however (cf. section 3 for more details)
the Galois closure
of $f$ yields another
quasi-\'etale map
$\phi : Z \ra X$. Since $\phi$ must be birational,
    it follows that $f$ is Galois and we have a contradiction to $Aut(X)
=  \{1\}$.

Is it possible to construct such a variety $X$ with  properties 1) - 4) ?

\end{itemize}
\end{oss}

The analysis of the Q.E.D. equivalence classes for Kuga-Shavel surfaces
is based on similar ideas, except that smooth ball quotients, or polydisk
quotients, have a residually finite fundamental group. The key result is then:

{\bf Theorem
\ref{quaternion} }
{\em Let $k$ be a real quadratic field, and let $\A$ be the
indefinite division quaternion algebra corresponding,
by Hasse's theorem, to a choice of $\SSS$
made as in \ref{S}.

Define $\FF_{\A}$ to be the family of subgroups
    $\Delta \subset \PP SL (2, \R ) \times  \PP SL (2, \R )$
   commensurable with a subgroup $\Gamma$
associated to a maximal order $ \RR \subset \A$.

Each $ \De \in \FF_{\A}$ acts freely
on $\HH^2$, and denote by
$S_{\De} : = \HH^2 / \De$
the corresponding algebraic surface.

Then the family of surfaces $ \{ S_{\De} | \De \in \FF_{\A} \}$
   is a union of Q.E.D. equivalence
classes.}

{\bf Question 4} Do there exist for each $ n \geq 2$ varieties obtained as ball
quotients, and which yield non standard varieties of general type ?

In the complex non K\"ahler world, things get complicated already
  in dimension $2$ for special surfaces.
In fact, a compact complex surface with odd first Betti number is non 
K\"ahler,
and in an appendix S\"onke Rollenske shows

{\bf Theorem  \ref{kodaira} }
{ \em Let $S$ be a minimal Kodaira surface. Then a
smooth surface
$S'$ is $Q.E.D.$ equivalent to
$S$ if and only if $S'$ is itself a Kodaira surface. Thus Kodaira surfaces
constitute a single $Q.E.D.$ equivalence class.}

Recall moreover that Kodaira dimension is known not to be deformation invariant
for compact complex manifolds which are not K\"ahlerian (due to some
examples orginating from
the work of Blanchard, cf. \cite{ue80}, \cite{ue82}
section 5. and also section 5 of \cite{ca02} for a more general description).

On the other hand, recently Claire Voisin (\cite{voi}, \cite{voi2})
has given a negative
answer to the so called Kodaira' s question whether a compact
K\"ahler manifold is always a deformation of a projective variety.

Her counterexamples however leave open the following more general question
(which may in turn also have a negative answer)

{\bf Question 5:}
Is a compact
K\"ahler manifold  always $\C$-Q.E.D. equivalent to an algebraic variety?

There is perhaps a reason why Q.E.D. equivalence maybe more meaningful
for algebraic varieties defined over $\C$.
    One should in fact keep in mind that every projective variety over $\C$ is
an algebraic deformation of
a projective variety defined over $\bar{\Q}$: this follows since
Hilbert schemes are defined
over $\Z$. Thus, the study of varieties defined over $\bar{\Q}$ could
play a key role
for the Q.E.D. problem (cf. the next sections for more questions).

A final observation is that also the classical questions of unirationality
can be seen through a different  perspective if we adopt Q.E.D. -equivalence:
for instance, the classical counterexamples to the L\"uroth problem are
(cf. remark  \ref{unir}) Q.E.D. -equivalent to projective space.

\bigskip

\section{The trivial equivalence relation  }

In this section we show, as already mentioned, the necessity,
in defining the Q.E.D. equivalence,
to put some restriction on the singularities of the varieties
that we consider.

\begin{df}
The t- equivalence relation
for Algebraic Varieties is the one
generated by

\begin{itemize}
\item
Birational equivalence
\item
Flat deformation with connected Base.
\end{itemize}
\end{df}

\begin{teo}
Two (irreducible) algebraic varieties are t-equivalent if and only if
they have the same dimension.
\end{teo}
\Proof
Let $Z^n$ be an irreducible algebraic variety of dimension $n$: then $Z^n$
is birational to a projective hypersurface $V^n_d \subset \PP^{n+1}$
(cf.\cite{hart}, I, prop. 4.9  ). In turn, by varying the coefficients of the
polynomial of degree $d$ defining $V^n_d$, we see that $V^n_d$
is deformation equivalent to the cone $ CW^{n-1}_d$ over a projective
variety $W^{n-1}_d \subset \PP^{n}$.
Since $ CW^{n-1}_d$ is birationally equivalent to
$\PP^1 \times W^{n-1}_d$, and we can easily show
that $X \sim_t X'$,$Y \sim_t Y'$ implies
$X \times Y \sim_t X' \times Y'$,
we infer by induction that our variety $Z$ is $t$-equivalent
to $\PP^{n-1} \times C^{1}_d$, where $C^{1}_d \subset \PP^2$
    is a plane curve of degree $d$. Obviously $C$ is deformation equivalent
to a rational nodal curve ${C'}^{1}_d \subset \PP^2$, which is birational
to $\PP^1$. Whence, $Z$ is $t$-equivalent
to $\PP^n$. Conversely, it is clear that $t$-equivalence respects
the dimension.

\QED

\begin{oss}
Actually, the proof holds more generally if we consider connected
algebraic varieties (i.e., reduced and pure dimensional).
\end{oss}

\section{Elementary properties of quasi-\'etale morphisms}

For the reader's benefit, recall that (cf. \cite{reid2})
a variety $X$ has canonical singularities iff

\begin{itemize}
\item
[1] $X$ is a normal variety of dimension $n$
\item
[2] $K_X$ is $\Q$-Cartier, i.e., there is a positive integer $r$ (the
minimal such
integer is called the {\bf index} of $X$)
such that
the Weil divisor $ r K_X$
    is Cartier. This means that the following holds: letting $ i : X_0
\ra X$ be the
inclusion of the nonsingular locus
of $X$, then  Zariski defined the canonical sheaf $K_X$ as $
i_{*}(\Omega^n_{X_0})$,
and we want that $ i_{*}(\Omega^n_{X_0})^{\otimes r}$ be invertible on $X$.
\item
[3] if $p : Z \ra X$ is a resolution of singularities of $X$,
$$ r K_Z = p^* ( r K_X)  + \sum_j a_j E_j ,$$ where the $E_j$'s are
the exceptional divisors,
and the $a_j$'s are {\bf non negative} integers ( $a_j \geq 0$).

\end{itemize}
It follows directly from the definition that, $\forall m \geq 0$, there is a natural
isomorphism between
$ H^0(X, m r K_X) : = H^0(X_0, (\Omega^n_{X_0})^{\otimes rm})$ and $
H^0(Z, m r K_Z)$.

Given $f : X \ra Y$ a quasi-\'etale morphism between varieties with
canonical singularities,
of degree $d$, let $r$ be a common multiple of the indices of $X,Y$.

By definition  there is an open set $X^0$ such that  $ X^0 \ra Y^0$ is \'etale,
and $ X - X^0$ has codimension $\geq 2$. W.L.O.G.  we shall assume always
that $ X^0 \subset X_0$, and we may further assume that  $ Y^0 \subset Y_0$.
Since $ X_0 - X^0$ and $ Y_0 - Y^0$ have both codimension $\geq 2$,
$ H^0(X, m r K_X)  = H^0(X^0, (\Omega^n_{X^0})^{\otimes rm})$ and the
same holds
for $Y^0$. Then $X$ and $Y$ have the same Kodaira dimension since
$$ H^0(Y, m r K_Y) \subset  H^0(X, m r K_X)  \subset H^0(Y, d m r K_Y) .$$

    \begin{oss}
Assume that $Y$ is smooth: then a quasi-\'etale morphism $f : X \ra
Y$ is \'etale.

\end{oss}

\Proof
In fact, $\pi_1(Y^0) \cong \pi_1(Y)$, thus there is an \'etale
covering $W \ra Y$
such that $X^0$ and $W^0$ are isomorphic. E.g. by Zariski's main
theorem, the birational
map induces a morphism $ f : X \ra W$. Moreover, cf. Theorem 7.17 of
\cite{hart},
and exercise 7.11 (c), page 171,  $X$ is the blow up of an ideal
sheaf $\I \subset \hol_W$
such that $ supp ( \hol_W / \I )$ has codimension $\geq 2$ and
contained in $W- W^0$.
But then, since the pull back of $\I$ is invertible, we contradict $
cod(X-X^0) \geq 2$
unless $\I$ is equal to $\hol_W$, that is, $X \cong W$.

\qed

\begin{prop}
Let  $f : X \ra Y$ be a quasi-\'etale morphism, and let $ f \circ g:
W \ra X \ra Y$ be the
Galois tower of $f$. I.e., over $\C$,  we let $ W^0 \ra X^0 \ra Y^0$
be the sequence of \'etale coverings
corresponding to the biggest normal subgroup contained in $\pi_1 (X^0
) \ra \pi_1 (Y^0)$,
and let  $ W \ra X$ be the corresponding finite normal ramified covering.

It is clear that  $ W \ra Y$ is quasi-\'etale, and we claim that if $X,Y$
have canonical singularities,
then also $W$ does.
\end{prop}

\Proof

$K_W$ is $\Q$-Cartier since $r K_W | W^0 = g^* ( r K_X | X^0)$ for
each positive
integer $r$, thus $r K_W  = g^* ( r K_X )$ as Weil divisors,
and it suffices to take $r$ such that $ r K_X$ is Cartier.

Let moreover $\pi_X : X' \ra X ,\pi_W : W' \ra W $ be  respective resolutions
such that $ g \circ \pi_W $ factors as $ \pi_X  \circ  g'$: then we have
$$ \pi_W^*( r K_W) = \pi_W^* g^* ( r K_X ) = (g')^* \pi_X^* (r K_X) = $$
$$=(g')^*
(r K_{X'} - \sum_j a_j E_j) = (r (K_{W'} -   R) - \sum_j a_j (g')^* E_j),$$
where $ a_j \geq 0$ and $R$ is the  ramification divisor (an
effective divisor) : whence
condition [3] is satisfied.

\qed

\begin{cor}\label{groupquotient}
Let  $f : X \ra Y$ be a quasi-\'etale morphism, with $X$ smooth and $Y$ normal.
    Then there is an
\'etale covering $W$ of $X$ and an action of a finite group $G$ on $W$,
    free in codimension $1$,
such that $Y$ is a birational image of $W/G$ under a small contraction.

Moreover, if $Y$ has canonical singularities, then also $W/G$ has canonical
singularities.
\end{cor}

\Proof
Let $ g : W \ra X$ be as in the previous proposition, thus $g$ is \'etale.

Letting $G$ be the Galois group, we obtain a birational morphism
$\phi : W/G \ra Y$
such that $ f \circ g = \phi \circ p$, $ p: W \ra W/G$ being the
quotient projection.

Now  the birational morphism $\phi$
induces an isomorphism $ W^0/G \ra Y^0$, i.e.,
    outside an algebraic set of codimension $\geq 2$,
which means that the contraction is small.

It is clear that $G$ acts freely in codimension $1$, whence it follows that
$K_{W/G}$ is $\Q$-Cartier, and since $\phi^*(K_Y) = K_{W/G}$ we obtain that
property [3] is satisfied.

\qed

The importance of the previous corollary lies in the fact that the
conditions that
$W$ be smooth, and $W/G$ have canonical singularities imposes restrictions
on the actions of the stabilizers $G_w$, for $ w \in W$.

Recall in fact that, by a well known Lemma by H. Cartan, if a finite
group $H$ acts
on a smooth germ $(\C^n, 0)$, then we may assume, up to a
biholomorphism of the germ,
that the action is linearized, i.e., $ H \subset GL(n,\C)$.

We have then (cf. \cite{reid2}) the following

\begin{oss}
Let $ H \subset GL(2,\C)$ be a finite group, and assume that

i) $H$ contains no pseudoreflections, i.e., $H$ acts freely in codimension $1$.

Then the germ $\C^2 / H$ has canonical singularities iff $ H \subset SL(2,\C)$.

A similar result does not hold true so simply for $n \geq 3$, cf. \cite{reid2},
exercise 1.10. page 352.
\end{oss}

Concerning the varieties with  Kodaira dimension $\K = -\infty$, it
is conjectured
that they are precisely the {\bf uniruled} varieties (cf.
e.g.\cite{kol} 1.12, page 189),
i.e., the varieties $X$ of dimension $n$ such that there exists a dominant
(and separable if char $\neq 0$) rational map $ Y \times \PP^1 \ra X$, where
$ dim Y = n-1$.

There is a related conjecture by Mumford  (see 3.8.1 of \cite{kol},
page 202)

{\bf Mumford's conjecture:}
    Let $X$ be a smooth projective variety of dimension $n$.
Then $X$ is separably rationally connected iff
$$H^0 ( (\Omega^1_X)^{\otimes m} )= 0  \  \forall m > 0.$$

In the context of the Q.E.D. problem, observe first that, by a result
of Fujiki and Levine
(cf. \cite{fuj}, \cite{lev}, and cf. also Chapter IV of \cite{kol}),
the class of
uniruled varieties is stable by deformation, at least over $\C$.

Then the following proposition (whose proof uses a precious
suggestion by Thomas Peternell)
ensures that the family of
complex uniruled varieties is stable by Q.E.D. equivalence

\begin{prop}\label{peano}

    Let $f : X \ra Y$ be a quasi-\'etale morphism, with $Y$ uniruled.
Then, if $X,Y$ have canonical singularities, also $X$ is uniruled.

\end{prop}

\Proof
By the theorems of Miyaoka and Mori, resp. Miyaoka (\cite{m-m} and
\cite{miya}, cf. also
\cite{kol}, Thm. 1.16, page 191), if we we have a smooth projective variety
$Z$ of dimension $n$, then $Z$ is separably uniruled
if and only if there is a covering family of  curves
$C_t, t \in T$, with $ C_t \cdot K_Z < 0$.

Let first $Z = \tilde{Y}$ be a resolution of $Y$.
    Then, if $Y$ is uniruled, there is such a family
of curves $C_t$ on $Z = \tilde{Y}$: let us push down this family to
a covering family of curves $D_t$ on $Y$.
By property 3. of canonical singularities, we obtain that
$ D_t \cdot K_Y < 0$.

Let $\Gamma_t$ be the family of proper transforms of the curves $D_t$
on $X$: since
$f$ is quasi- \'etale, it follows from the projection formula that
$\Gamma_t \cdot K_X < 0$.

Let us consider a resolution $X'$ of $X$: then by \cite{bdpp} it suffices to
show  that $K_{X'}$ is not pseudoeffective.

Else, if $L$ is an ample divisor on $X$, then for $ m, N >> 0$
the linear system  $ | N (m  K_{X'} +  L )|$ is effective and  big (yields
a birational map). In particular, if we denote by $\Delta_t$ the proper
transform of a general $\Gamma_t$, the intersection number
$ (m  K_{X'} +  L )  \cdot \Delta_t \geq 0$.

This is however a contradiction, since for $m$ sufficiently divisible
and by the
projection formula it follows that $ | N (m  K_{X'} +  L )| = \pi^* |
N ( mK_X + L)|$.
Whence, by the projection formula,  $ (m  K_{X'} +  L )  \cdot \Delta_t
= ( mK_X + L) \cdot \Gamma_t$. This last number is however negative for
$ m >> 0 $, a contradiction.

\qed

It is not clear  whether the following question should have a
positive answer:

{\bf Question 6 :} If a variety  $X$ (we assume canonical singularities throughout) 
has   Kodaira dimension $\K = -\infty$, then  is there a quasi-\'etale morphism $f : Z
\ra X'$ where $Z$ is ruled and $X'$ is a deformation of $X$?

\begin{oss}\label{unir}
Observe that a general cubic threefold $X \subset \PP^4$ is unirational
but not rational, and that a smooth quartic threefold $Y \subset \PP^4$
is unirational but not rational. By Lefschetz' theorem they are both simply
connected (as Koll\'ar rightly reminded  me), whence they have no
nontrivial quasi-\'etale cover. However, we can deform $X$ to a cubic threefold
$X'$ with a double point, respectively $Y$ to $Y'$ with a triple point.
Both $X', Y'$ have canonical singularities and are rational,
whence a positive answer
to the above question in this special case.

In the case of the quartic threefold, by \cite{I-M}, any smooth such manifold
is not rational; this class of smooth Fano manifolds is stable by
deformation only if we restrict ourselves to the condition that the fibres be smooth
and projective.

\end{oss}

Already for conic bundles it is not clear whether question 6 has a
positive answer.

The second remark above shows
  that the following
stronger question has a negative answer:
if a variety $X$ has   Kodaira dimension $\K = -\infty$,
   is there a quasi-\'etale morphism $f : Z \ra X'$ where $Z$ is ruled,
$Z$ and $X'$ are smooth,
and $X'$ is a deformation of $X$?

This question is somehow related to a
stronger version of the Mumford conjecture :

{\bf Quasi-\'etale unirationality questions:}  Let $X$ be a smooth
projective variety of dimension $n$.
Then
$$H^0 ( (\Omega^1_X)^{\otimes m} )= 0  \  \forall m > 0$$
implies that $X$ is quasi-\'etale equivalent to $\PP^n$ and unirational?

Is the condition $$H^0 ( (\Omega^1_X)^{\otimes m} )= 0  \  \forall m > 0$$
invariant by deformation of smooth projective varieties ?

\begin{oss}
It is easy to see that the property $H^0 ( (\Omega^1_X)^{\otimes m} )= 0  \
\forall m > 0$ is invariant for quasi-\'etale morphisms between
smooth varieties.

It would be interesting to extend this condition for varieties with canonical
singularities.
\end{oss}

\begin{oss}s
Koll\'ar constructed (see e.g. Chapter V, section 5 of \cite{kol},page 273 and
foll.) examples of complex Fano varieties (these are
rationally connected) which are not ruled. Are these counterexamples
to question 6 ?
\end{oss}

\section{Proof of the main theorems \ref{kaehler} and \ref{algebraic}}

\Proof of THEOREM \ref{kaehler}

We must show that, if $S$ and $S'$ are smooth K\"ahler surfaces which
have the same
Kodaira dimension $\K \leq 1$,   then $S$ and
$S'$ are
$\C - Q.E.D.$-equivalent.

We proceed distinguishing the several cases,
according to the value of the  Kodaira dimension
$\K$.

{\underline {$ \K $$ = - \infty$.}}

In this case $S$ is projective algebraic (since $p_g(S)=0$ implies,
by the K\"ahler assumption,
that $b^+ > 0$, whence there is a positive line bundle))
and it suffices to show that $S$ is A.Q.E.D. - equivalent to $\PP^1
\times \PP^1$.

But $S$ is birational to a product $C' \times \PP^1$, and the curve
$C'$ is deformation
equivalent to a hyperelliptic curve $C$. Let $\iota$ be the
hyperelliptic involution,
and let $j: \PP^1 \ra \PP^1$ the involution such that $j(z) = -z$.

We have an action of $\Z/2$ on $C \times \PP^1$
provided by $\iota \times j$, which has only isolated fixed points. Set
$ X : = (C \times \PP^1)/\Z/2 $: there is a fibration $f : X \ra C
/(\Z/2) \cong \PP^1$,
thus by Noether's theorem $X$ is birational to $\PP^1 \times \PP^1$.

\QED

\begin{oss}\label{productinv}
More generally, let $C_1, C_2$ be hyperelliptic curves, so that we
have an action
of $(\Z/2)^2$ on $C_1 \times C_2$, and let us consider the diagonal embedding
of $(\Z/2) \subset (\Z/2)^2$. Set  $X:= C_1 \times C_2$,  $Y:= (C_1
\times C_2)/ (\Z/2)$
and observe that $(C_1 \times C_2)/(\Z/2)^2 \cong (\PP^1 \times \PP^1)$.

We have $ f: X \ra Y, p : Y \ra \PP^1 \times \PP^1$, and $Y$ has only nodes
as singularities, while $f$ is quasi-et\'ale, so that $X$ is A.Q.E.D.
equivalent to $Y$.
On the other hand, the branch locus $B'$ consists (if we denote by
$g_i$ the genus of
$C_i$) of the union of $(2g_1 + 2)$ vertical lines with
$(2g_2 + 2)$ horizontal lines.

Let $S$ be a double cover of $(\PP^1 \times \PP^1)$ branched on a smooth
curve $B$ of bidegree $(2g_1 + 2, 2g_2 + 2)$: since $B$ is a
deformation of $B'$,
$S$  is a deformation of $Y$ hence it is  A.Q.E.D. equivalent to $X$.

Observe that the composition of the double cover $S \ra (\PP^1 \times \PP^1)$
with the second fibration yields a fibration of hyperelliptic curves
of genus $g_1$.

In the particular case where $g_1=g_2=1$ the above argument
    shows that a product of elliptic
curves is  A.Q.E.D. equivalent to a K3 surface which is a double cover of
$(\PP^1 \times \PP^1)$.
\end{oss}

We proceed now with the proof of the next cases:

{\underline {$ \K $$ = 0$.}}

Recall that a minimal compact complex surface with Kodaira dimension
$0$ is either

a) A complex torus

b) A K3 surface

c) An Enriques surface

d) A hyperelliptic surface

e) A Kodaira surface.

In case e), $b_1(S) =3$ if $S$ is primary, else $b_1(S) =1$ if $S$ is
secondary. In particular, $S$ is not Ka\"hler.

In cases c) and d) $S$ is projective and it has a finite unramified covering
$f:  S' \ra S$ where:

\begin{itemize}
\item
c) if $S$ is Enriques, $deg(f)=2$ and $S'$ is a K3 surface
\item
d) If $S$ is hyperelliptic, $deg(f)| 12$ and $S'$ is a product of
elliptic curves.
\end{itemize}

By virtue of the previous remark \ref{productinv} the proof is concluded
since it is well known that all complex tori are deformation equivalent,
and Kodaira proved (\cite{kod1}) that every K3 surface is deformation
equivalent
to a smooth quartic surface in $\PP^3$.

\QED

{\underline {$ \K $$ = 1$.}}

Recall that a complex surface of Kodaira dimension $ 1$ is properly
(canonically) elliptic,
i.e., it admits a (pluri-) canonical elliptic fibration $f : S \ra B$.

Since we are interested in the case where $S$ is K\"ahler we can disregard
the case where $b_1(S)$ is odd (and $f$ is then an elliptic quasi-bundle).

Step I) We show first, replacing $S$ by a finite unramified covering, that we
may assume that $S$ has
an elliptic fibration without multiple fibres, unless we are in the following

{\bf EXCEPTIONAL CASE* :}  $f : S \ra \PP^1$ has at most two multiple fibres,
    with coprime multiplicities, and $S$ is simply connected.

PROOF of I: we use (cf. \cite{cko}, lemma 3 and Theorem A for a
similar idea, and also
    e.g. \cite{cat03} 4.1, 4.2 ) the orbifold fundamental group sequence
$$ \pi_1 (F) \cong \Z^2 \ra \pi_1(S) \ra \pi_1^{orb}(f) \ra 1 $$
where $F$ is a smooth fibre $F_b:= f^{-1}(b)$, $F_{b_1}, \dots F_{b_r}$ are the
multiple fibres, of respective multiplicities $m_1, \dots m_r$, and
$\pi_1^{orb}(f)$
is defined as the quotient  $\pi_1 ( B - \{b_1, \dots b_r\})/ <<
\gamma_1^{m_1}, \dots
\gamma_r^{m_r}>>$ of  $\pi_1 ( B - \{b_1, \dots b_r\})$ by the
subgroup normally generated
by the respective $m_j$-th powers of simple geometric loops
$\gamma_j$ around the respective points $b_j$.

Note that the image of $\gamma_j$ inside  $\pi_1^{orb}(f)$ has order
precisely $m_j$
unless we are in the exceptional case where $B \cong \PP^1$ and $ r \leq 2$
    (for $r=1$ the group is trivial, else
for $r=2$ it is cyclic of order $= G.C.D. (m_1, m_2)$).

It is also known (cf. \cite{cko}, loc. cit.) that, if we are not in
the exceptional case,
there is a finite quotient $G$ of $\pi_1^{orb}(f)$ where the image of
each  $\gamma_j$
has order precisely $m_j$. To this surjection corresponds an
unramified covering
such that the normalization of the pull back of $f$ is an elliptic fibration
without multiple fibres.

In the exceptional case with two multiple fibres
we may  take a cyclic cover $\PP^1 \ra \PP^1$
of order $= G.C.D. (m_1, m_2)$, branched on the two  points
corresponding to the
multiple fibres, so that
the normalization $S'$ of the pull-back has two multiple fibres whose
multiplicities
$m'_1, m'_2$ are coprime, whence $\pi_1^{orb}(f')$
is trivial, and $\Z^2 \ra \pi_1(S')$ is surjective.
Since we assume the first Betti number to be even, and since by \cite{dolg},
Theorem on page 137,
$\pi_1^{orb}(f)$ contributes here to the torsion subgroup of $H^1(S,
\Z)$, we infer
that also $S'$ has even $b_1(S')$. But if $b_1(S') = 2$ $S'$ is a
trivial family,
contradicting $Kod(S)=1$. Thus $\pi_1(S')$ is finite Abelian and passing to
the universal cover we find ourselves in the EXCEPTIONAL CASE*
(actually, one can indeed
show that $S'$ is itself simply connected).

Step II) Let $q$ be the irregularity of $S$: if $ q =  g(B) + 1$,
where $g(B)$ is the genus
of $B$, then $S$ is a product $B \times F$, where $F$ is a smooth
fibre, and $B$ is
a curve of genus $g \geq 2$. Thus we have only one
A.Q.E.D.-equivalence class, to which we
shall show that all the other cases are $\C$-Q.E.D. equivalent.

Step III) Consider now any two numbers $g,q$ with $ g \geq 2q$:
    then there is a ramified
double covering $ j : B \ra C$ where $B$ has genus $g$ and $C$ has genus $q$ .
Consider then a product $B \times F$ as before, and use the trick of
\ref{productinv},
to obtain an elliptic surface without multiple fibres
$ S \ra C$ with $q(S) = q$, $p_g(S) = g-q \geq q$.

Step IV) Assume that $S$ is an elliptic surface without multiple fibres
and with topological
Euler number $e(S)=0$: then, by the Zeuthen-Segre theorem
(cf. \cite{bpv},11.5, page 97), all the fibres of $ S \ra C$ are smooth.
It follows then that the $j$-invariant is constant, since $ j: C \ra
\C$ is holomorphic,
thus all the fibres are isomorphic and we have a holomorphic fibre bundle.

In this case the Jacobian elliptic surface $J$ associated to $S$ has
an \'etale cover which is a product (cf. e.g. \cite{bpv}, (2) page 143,
since there is an \'etale cover of $C$ which pulls back a principal bundle
with a section, cf. also \cite{f-m}, section 1.5.4)
and by theorems 11.9 and 11.10 of \cite{KodIII} (cf. also Theorem
11.1 of \cite{bpv})
it follows that,
since $b_1(S)$ is even, then $S$ is a complex deformation of $J$
and thus $S$ is in the $\C-Q.E.D.$ class II).

Step V) We may then assume that, if the  elliptic surface $S$ has no
multiple fibres,
then it has topological
Euler number $e(S) > 0$: by the Noether formula $ 12 \chi(\hol_S) =
K^2_S + e(S)$,
and since $K^2_S = 0$,
this means that $ p_g(S)  \geq q(S)$.

We use now theorem 7.6 of \cite{f-m}, asserting that two complex
elliptic surfaces without
multiple fibres, with $e(S) > 0$, and with the same $ q, p_g$ are
complex deformation equivalent.

By Step III) follows then that any such is $\C-Q.E.D.$-equivalent to
a product $B \times F$,
and we have therefore shown that there is only one $\C-Q.E.D.$-class,
unless possibly if we are in the EXCEPTIONAL CASE*,
which we treat next.

Step VI) Assume now that $S$ is simply connected, and that $f : S \ra
\PP^1$ has
multiple fibres, and at most two, of coprime orders $ 1 \leq m_1 < m_2$.
A further invariant of $f$ is the geometric genus $p_g(S) =
\frac{1}{12} e(S) -1.$

Two such surfaces with the same  invariants  $p_g(S), m_1, m_2$ are known
to be complex deformation equivalent (cf. \cite{f-m}, Theorem 7.6).

Therefore it suffices to find, for each choice of $p_g(S), m_1, m_2$ as above,
one such exceptional elliptic fibration $f: S \ra \PP^1$ which is
QED-equivalent to
one without multiple fibres.

To this purpose it suffices to find a divisor $D$ contained in a finite union
    of fibres and which is a disjoint union of (connected)
(-2)-curves configurations $D_1, \dots
D_k$,  such that the open surface $S^0 := S- D$ has now at least
three multiple fibres.
Because in this case there exists an unramified covering ${S'}^0$ of
$S^0$ which yields an
elliptic fibration without multiple fibres (again by the method of considering
the orbifold fundamental group exact sequence for the open surface,
as in \cite{cko}).
And $S^0$ is the complement of a finite set in the surface $X$ with
rational double points
obtained by contracting the  (-2)-curves configurations $D_1, \dots
D_k$ to points.
Likewise,  ${S'}^0$ is the complement of a finite set in a surface $Y$
with rational double points, mapping in a quasi-\'etale way to $X$,
thus $S$ is $\C-QED$-equivalent to the minimal resolution $S'$ of $Y$,
which has an elliptic fibration without multiple fibres.

We are then reduced to show the existence, for given $p_g(S), m_1, m_2$, of
an exceptional elliptic surface with those invariants and moreover with two
singular fibres (only one suffices  indeed if $m_1>1$) whose extended
Dynkin diagram
is of type $ \tilde{D}_n , n \geq 4$. or of type $\tilde{E}_n , n = 6,7,8.$

Observe that the value of $p_g(S)$ is determined by $e(S)$, and that
logarithmic
transformations do not change $e(S)$.

We are then reduced to showing the existence of a simply connected
elliptic  fibration
(over $\PP^1$) with at least two
singular fibres  whose extended Dynkin diagram
is of type $\tilde{D}_n , n \geq 4$, or of type $\tilde{E}_n , n = 6,7,8.$

As in remark \ref{productinv} lets us consider a double cover of
$\PP^1 \times \PP^1$ branched on a divisor $B$ of bidegree $(2 g + 2, 4)$.
If $B$ is smooth we get a simply connected elliptic surface $S$ with
$p_g(S) = g$. If $g \geq 1$ it is easy to show that we may obtain a
branch curve
    $B'$ with two ordinary triple points:
then the double covering surface $S'$ gets two singular fibres of type
$\tilde{D}_4 $.

If $g=0$ we obtain $B'$ as the union of four divisors: $L_1, L_2$ of bidegree
$(0,1)$ and $D_1, D_2$ of bidegree $(1,1)$.

Viewing in fact $\PP^1 \times \PP^1$ as a smooth quadric in $\PP^3$, letting
    $L_1, L_2$  be two disjoint sections, and fixing $P_1 \in L_1, P_2
\in L_2$, then
we choose $D_1, D_2$ as two general plane conic sections through $P_1, P_2$.

\QED for Theorem \ref{kaehler}

PROOF of Theorem \ref{algebraic}.

We essentially rerun the proof of \ref{kaehler}, mutatis mutandis.

For $\K = - \infty$ the K\"ahler surfaces are algebraic and the proof
is already there.

For $\K = 0$ we simply have to observe that

1) any abelian surface is an algebraic
    deformation of a product of elliptic curves

2) two algebraic K3 surfaces are algebraic deformation of each other.

Statement 1) is easy, in any dimension $n$, since Abelian varieties with
a polarization
of type $(d_1, d_2, \dots d_n)$ are parametrized by a quotient of
the Siegel upper halfspace, and a product of elliptic curves clearly admits
such a polarization.

     The case of K3 surfaces is  similar and requires the Torelli
theorem (\cite{ps-s},
cf. also \cite{K3}, expose' XIII): again
we have an irreducible subvariety parametrizing the projective K3
surfaces with a
    pseudopolarization of degree $d$, and inside this
family we find the special Kummer surfaces, i.e., more precisely, the
K3 surfaces
obtained as the minimal resolution of
the Kummer surface of a product of elliptic curves.

Let's consider now the case $\K = 1$.

Steps I, II, III are identical.

Step IV: assume  $e(S) = 0$: then $ f : S \ra C$ is a holomorphic bundle
and there is ( cf. \cite{bpv}, page 143)
an \'etale covering of the base $C$ such that the pull back
is a principal holomorphic bundle with cocycle $\xi$ whose
cohomological invariant $ c(\xi) = 0$ (else, by \cite{KodI} theorem 11.9,
cf. also \cite{bpv} Prop. 5.3,
page 145, $ b_1(S) \equiv 1 ( mod 2)$, contradicting the algebraicity of  $S$).

Let $F$ be the fibre of $ f$: then (cf. also \cite{f-m} page 92)
there is a finite
homomorphism $ \pi_1(C) \ra F$ classifying $f$, and taking the
associate \'etale
cover $C'$, we obtain an \'etale covering $ S' \ra S$ which is indeed
a product.

Step V : recall  that a Jacobian
    elliptic surface
is algebraic.  As shown by Seiler  (cf. \cite{sei}, and also
 \cite{kas} or \cite{mir} for an  introduction to the subject)
all  Jacobian elliptic fibrations which are not a product and have
the same invariants $q(S), p_g(S)$ belong to an irreducible algebraic family.

Therefore, any Jacobian elliptic surface is an algebraic deformation
of some Jacobian surface with constant invariant $j$ obtained from construction
4.1 as in step III).

Let us use the fact that the base  space of a maximal family of
algebraic elliptic surfaces is a finite
covering of the corresponding base space of the corresponding family of
Jacobian elliptic surfaces (cf. \cite{sei}, and also \cite{f-m} prop.
5.30, page 93).
This is derived from Kodaira's theorem 11.5 of \cite{KodI} asserting that if
$S$ is an algebraic elliptic fibration without multiple fibres, then
the corresponding cohomology class $\eta$ is torsion, and conversely.

We conclude that an algebraic elliptic fibration without multiple fibres
$f : S \ra C$ is algebraic deformation of one with constant moduli
and with multiplication by $ \pm 1$. Whence, a double  etale covering
of the base yields a  double  etale covering of $S$ which is a
holomorphic bundle.
We are done by Step IV.

    Step VI): we are in the exceptional case where
$ f : S \ra \PP^1$ has multiple fibres, of coprime multiplicities $ m_1 , m_2 $
with $ 1 \leq m_1 < m_2 $.

We are done once we can show the validity of the following

{\bf Claim:} Assume that we have an algebraic exceptional elliptic
surface $S \ra \PP^1$, i.e., with $ r \in \{1,2\}$ multiple fibres,
of coprime multiplicities $ m_1 , m_2 $
with $ 1 \leq m_1 < m_2 $. Then there exists an algebraic deformation of $S$
yielding a surface $S'$ with $3-r$ singular fibres of
type $\tilde{D}_n , n \geq 4$ or  $\tilde{E}_n , n = 6,7,8.$

{\em Proof of the claim.}
We argue as in Step VI of the previous Theorem \ref{kaehler},
using the  characterization of algebraicity of logarithmic transforms
given in \cite{f-m}, Lemma 6.13, page 106, and which is a translation
in complex geometry of the theory of Ogg-Shafarevich (cf. \cite{shaf1}, \cite{dolg}).

Let us then consider an algebraic exceptional elliptic  surface $\phi : S
\ra \PP^1$, and let
$\psi: X \ra \PP^1$, be its Jacobian
fibration. By doing the inverse of logarithmic transformations, we
obtain an  elliptic fibration
$f : Y \ra \PP^1$, whose Jacobian fibration is $X$.

The following is the content of the cited Lemma 6.13 of \cite{f-m}.

\begin{oss}
Let $f' : Y' \ra B$ be another  elliptic surface whose Jacobian 
fibration is $X$,
and let
$\phi ': S' \ra B$ be obtained from $Y'$ via the same logarithmic
transformations as the ones constructing $S$ from $Y$
(i.e., at the  fibres over the same points, and with the same
associated torsion bundles): then $S'$ is algebraic iff and only if
the difference of the corresponding elements in the
classifying group $ H^1(\PP^1, \hol (X^*))$ ($\hol (X^*)$ is Kodaira's
sheaf of groups of  local holomorphic sections) is torsion.
\end{oss}

Argueing as in Step VI of \ref{kaehler}, there is an algebraic
1-parameter family of
Jacobian elliptic surfaces $X_t, t \in T$ containing the given $X$
and a special one $X_0$ which
has $3-r$ singular fibres of type $\tilde{D}_4$.

Let us treat now the case $ 3-r=1$.
Let us consider the  elliptic surfaces $Y_{t, w}, w \in W_t,$ without 
multiple fibres
having some $X_t$ as Jacobian elliptic surface and the
  family of
single logarithmic transforms of the
surfaces $Y_{t,w}$ : this family is parametrized by a complex variety
$Z$ where $ Z
\ra W$ has  pure and irreducible 1-dimensional fibres.

Inside this family we consider the subfamily of the algebraic 
elliptic surfaces: these form
a countable union of subvarieties
fibred over our irreducible curve
$T$, whence, up to replacing
$T$ by an irreducible  finite covering of it, we may find,
  given our initial $S$,  a 1-parameter complex family
$S_t$ containing $S$ and
such that the corresponding family of Jacobian surfaces is $X_t$.

This shows that we have such a complex family of algebraic surfaces.
In order to show that we have an algebraic family we only need to observe
that our elliptic surfaces have all a multisection $D$ of a fixed degree,
whence for very large $n, m  \in \N$ $ |n D + m F|$ is very ample on
each surface, and we get a non trivial complex curve in a Hilbert scheme of
projective surfaces. We only  need to remark that if two points of a Hilbert scheme
are joined by a complex curve, they are also joined by an algebraic curve.

The argument for the case $ 3-r=2$ is entirely similar, whence our
claim is proven,
together with Theorem \ref{algebraic}.

\QED

\section{Remarks and questions.}

We want in this section to add more questions, and some comments  regarding
some of the questions previously posed in the introduction.

1) Reduction to varieties over $\overline{\Q}$ can be thought of as the
distinguishing feature
between algebraic and K\"ahler varieties. In fact, if we have a
smooth projective
variety $X \subset \PP^n_{\C}$, we get a corresponding point $[X]$ of
a Hilbert Scheme
$\HHH$. Since $\HHH$ is defined as a closed algebraic set
in an appropriate Grassmannian $G$ by rank equations of certain
multiplication by monomials,
$\HHH$ is defined over $\Z$, and an irreducible component containing $[X]$
contains a dense set of points defined over $\overline{\Q}$, whence we obtain a
smooth projective
variety $Y \subset \PP^n_{\C}$, where $Y$ is defined over
$\overline{\Q}$ and is an
algebraic deformation of $X$. Assume that $Y$ is defined over a number field
$K$: then the theory developed so far suggest to consider the quasi-\'etale
generalization of Grothendieck's fundamental group, which should play
an important
role (in case of canonical models of varieties of general type, i.e.,
of varieties $X$ with canonical
singularities  with $K_X$ ample, this is exactly the Grothendieck
fundamental group of the smooth
locus of $X$).

It would be also interesting to enlarge our equivalence relation as
to include, for varieties defined over a number field,
also the action of the absolute Galois group
(thus for instance considering a variety  $ / \C$ and its
complex conjuate as equivalent).

{\bf Question 7:}
  For which  classes of algebraic varieties  is A.Q.E.D.-equivalence the
same as the weaker
$\C$-Q.E.D. equivalence ?

2 ) In the case of uniruled varieties the Q.E.D. question is strictly
related to  the question of
"generic" splitting of normal bundles for the curves of a covering
family of rational curves
(here, "generic" stands not only for the generic curve of the family,
but also for a general
deformation of the given variety).

4) what is the t-equivalence of compact complex manifolds? (this is
hard since for instance we do not know
all the compact complex surfaces).

{\bf Question 8} Assume that $S =  B^2 / \Gamma$ is a compact minimal
smooth surface which is a ball
quotient, (equivalently, by Yau and Miyaoka's theorem , cf.
\cite{yau} and \cite{miy}, $ K^2_S = 9 \chi
(\hol_S)$). Does there exist, as in the case of Kuga-Shavel surfaces,
  such a group
$\Gamma$ such that every group $\Gamma ''$ commensurable
with $\Gamma$ is either torsion free (it acts freely), or it has a
fixed point $z$ where
the (finite) stabilizer $\Gamma ''_z$ has a tangent representation
not contained in
$ SL ( 2, \Z)$ ?

Fritz Grunewald suggested that there such examples should indeed exist,
more precisely that there are such groups $\Gamma$ such that every
group $\Gamma ''$ commensurable
with $\Gamma$ equals $\Gamma$, and such examples
should be found among the
ones of Deligne-Mostow (cf. \cite{D-M}).

{\bf Question 9}(Lucia Caporaso): which are the Q.E.D. equivalence class of
Kodaira fibrations ?

\section{ Appendix due to Fritz Grunewald: \\
Q.E.D. classes constructed from quaternion algebras. }

As in \cite{shav} (cf. also \cite{shim}, chapter 9) we consider a
division quaternion algebra
$\A$ with centre a totally real number field $k$. For simplicity, we may
further assume $k$ to be a  real quadratic field.

We assume further that $\A$ is totally indefinite, which  means that for each
of the two embeddings
$ j : k \ra \R$ $\A$ does not ramify, i.e.,  $\A \otimes_j \R  \cong
M( 2, \R) :=  Mat ( 2 \times 2, \R).$

As usual, denoting by $\hol_k$ the ring of integers of $k$, and by
$k_{\PPP}$ the local field which is the completion
of the localization $\hol_{\PPP}$ of  the ring
$\hol_k$ at a prime ideal $ \PPP$,
one
considers the set of primes where
$\A$ ramifies, i.e., the subset

$$ \SSS (\A) : = \{ \PPP  \in Spec (\hol_k) | \A \otimes_k k_{\PPP}
{\rm is \ a \ skew\  field \ } \} .$$

By the classical results of Hasse (which are exposed for instance in the
book \cite{weil}, cf. especially Th. 2 of Chapter XI-2, and Theorem 4,
section 6 of chapter XIII ) we know that (cf.  also \cite{shim}, section 9.2,
pages 243-246)

\begin{itemize}
\item
s-1) The cardinality of $ \SSS (\A)$ is finite and even (and nonzero since
$\A$ is a division algebra)
\item
s-2) $\A$ is completely determined by  its centre $k$ and by $ \SSS (\A)$
\item
s-3) for each choice of  $k$ and of a set $ \SSS \subset Spec (\hol_k)$
with even cardinality, there is a quaternion algebra over $k$
with $ \SSS (\A) = \SSS .$

\end{itemize}

\begin{oss}
Usually one considers inside $ \SSS (\A)$ also the places at infinity
(embeddings of $k$ into $\R$), and the result holds more generally.
Since however we assume the quaternion algebra to be totally indefinite,
there are no ramified places at infinity.
\end{oss}

Let now $\RR \subset \A$ be a maximal order
(an order, cf. \cite{weil}, def. 2 page 81, is a subring which is
a $\Q$-lattice for $\A$)  and  consider the
group
$$ \Gamma(1) : = \{ a \in \RR | \  nr (a) = 1 \}, $$
where $nr$ denotes the reduced norm (\cite{weil}, IX-2).

The following facts are also well known (cf. \cite{shav}, section 1,
and \cite{shim}, 9.2)

\begin{oss}
1) If $k$ is a quadratic field, and $j_1, j_2$ are the two embeddings
   $ k \ra \R$, then $ \Gamma(1) \cong
(j_1 \times j_2) (\Gamma(1) ) \subset  SL (2, \R ) \times  SL (2, \R ).$

2) The image $\Gamma$ of $ \Gamma(1)$ in
$ \PP SL (2, \R ) \times  \PP SL (2, \R )$ is isomorphic to
   $ \Gamma(1) / \{\pm 1\}$.

3) $\Gamma$ operates properly discontinuously with compact
projective quotient
on the product $ \HH^2$ of two upper -halfplanes ( $\HH :=
\{ z \in \C |  Im \ z > 0 \}$).

4) the action of $\Gamma$ on  $ \HH^2$ is irreducible; whence,
if the action of $\Gamma$ is free, then the projective surface
$ X : =  \HH^2 / \Gamma$ is strongly rigid (cf. \cite{j-y},  and also
\cite{cat0}),
i.e., every surface $S$ with the same Euler number as $X$
and with isomorphic
fundamental group $\pi_1(S) \cong \pi_1 (X)$ is either
biholomorphic to $X$ or to the complex conjugate
surface $\overline{X}$.

5) Assume the  quotient $ X : =  \HH^2 / \Gamma$ to be smooth:
then its first Betti number equals zero (proposition 2.1 of \cite{shav},
which follows by the theorem of Matsushima and Shimura), and
by Hirzebruch's proportionality principle, we have
$ e (X) =  2 + b_2(X) = 4 ( 1 + p_g(X))$.
\end{oss}

\begin{lem}\label{span}
   Let $ \Gamma " \subset \Ga (1)$ be a finite index subgroup.
Then the $\Q$-linear span
of $\Ga ''$ equals
$\A$ .
\end{lem}

\Proof
Replacing $\Ga ''$ by a subgroup of finite index (since $\Ga (1)$
is finitely generated), we may assume that
$\Ga ''$ is a normal subgroup of  $\Ga (1)$ and invariant
by the involution of $\A$ sending an element to its conjugate
(it has $k$ as set of fixed points).

As a first step, let's now prove that

i) the $\Q$-linear span $G$
of $\Ga ''$ contains $k = k \cdot  1$.

In order to show this, let us consider an element $\alpha\in k$
subject to the conditions:
\begin{itemize}
\item $\alpha$ is totally positive,
\item $k ( \sqrt \alpha )\otimes_k k_{\PPP}$ is a field for every
$\PPP\in \SSS (\A)$,
\item there is a finite place ${\mathcal Q}\notin \SSS (\A)$ of $k$
such that $k_{\mathcal Q}$ is an extension field of degree $2$ over the
corresponding completion of $\Q$ and
such that $\alpha$ is not a square in $k_{\mathcal Q}$,
\end{itemize}
and set
$$L:=k(y)\qquad \qquad (y^2=\alpha).$$
The existence of such an $\alpha$ is guaranteed by the weak approximation
theorem (cf.  18.3 exercise  2 , page 351 of \cite{pie}) which says that $k$ is
dense when diagonally embedded into the direct
product of any finite  subset of the set  of its completions. Note 
also that the set
of non-squares is open in any completion of $k$.

The field $L$ has the following properties:
\begin{itemize}
\item $L$ is totally real,
\item $L$ is isomorphic to a subfield of $\A$,
\item the only subfield of $L$ which is of degree $2$ over $\Q$ is $k$.
\end{itemize}
The second condition follows from the second condition on $\alpha$,
cf. prop. 4.5 of \cite{shav}.

Then
we have a chain of degree 2 extensions $ \A \supset  L \supset k$,
and $k$ is the only
quadratic subfield of $L$, in fact the Galois group of
the splitting field of $L$ over $\Q$  is either cyclic
of order 4 or is the Dihedral group $D_4$.

Now, consider $\Ga '' \cap L$: we claim that
$\Ga '' \cap L$ contains a nontrivial
infinite cyclic subgroup generated by a unit $\epsilon$.

In fact, the maximal order $\RR$ of $\A$ intersects $L$
in an order $\BB$ of $L$, and by Dirichlet's  Theorem (cf.
e.g. Theorem 5 of section II, 4 of
\cite{b-s}, page 113) the group of units in  $\BB$ has rank 2;
since by the same theorem the group of units in $\BB \cap k$
has rank 1, the group of units of norm 1 in $\BB$ contains an infinite
cyclic group, which in turn intersects $\Ga ''$ in an infinite cyclic group.

Note that $\epsilon \notin k$, since for elements of $k$ the norm is just
given by the square, and $ \epsilon \neq \pm 1$.

Thus $\epsilon ,
\overline{\epsilon} \notin k $
but, clearly,  $ (\epsilon +  \overline{\epsilon} ) \in G \cap k $
and  we claim that $ (\epsilon +  \overline{\epsilon} ) \notin \Q $.

Otherwise, if  $ (\epsilon +  \overline{\epsilon} ) \in \Q $,
since $ \epsilon \cdot   \overline{\epsilon} = 1$
it follows that  $\epsilon ,  \overline{\epsilon}$ belong to a 
quadratic extension of
$\Q$. But this quadratic extension, being contained in $L$, would
then equal $k$, contradicting our previous assertion.

We  conclude then that $ (\epsilon +  \overline{\epsilon} )
  \in k \setminus \Q $, thus $G$ contains $k $.

ii) Now that we know that the $\Q$-linear span $G$
of $\Ga ''$ contains $k \cdot  1$, we can show that $G$ is a field.
In fact we observe  that $G$
is a ring, which is invariant by the involution of $\A$: thus
if $G$ contains $x$, it contains also $ x ^{-1} = \bar{x} \cdot rn(x)^{-1} $.

iii) Set $d: =  dim_k (G)$: then $ d | 4$, and if $ d=4$ there is
nothing to prove. If instead $ d \leq 2$, then $G$ is commutative,
hence also $\Ga ''$ is commutative.

  This gives however two contradictions:

1)  since we know that there exists a  finite index subgroup of
  $\Ga ''$ which is infinite,  operates freely on $\HH^2$, and has a 
quotient $X$
having a finite homology group.

2)  since we know
that $\Ga ''$ is Zariski dense in $ SL (2, \R ) \times  SL (2, \R )$,
by \ref{dense}.

\QED

\begin{lem}\label{dense}
Any finite index subgroup $\Ga ''$ of $\Ga(1)$ is
Zariski dense in $ SL (2, \R ) \times  SL (2, \R )$.

  $ \A^* : = \A \setminus \{0\}$ is Zariski dense in $ GL (2, \R ) 
\times  GL (2, \R )$.
\end{lem}

\Proof
The first assertion is a special case of a general theorem by Armand Borel
(cf. \cite{bor}).
In fact $\Ga ''$ has the Selberg property since the quotient
$\HH^2 / \Ga''$ is compact, hence $\Ga ''$ is not contained
in any proper subgroup having a finite number of components.

The second assertion follows immediately since  $\Ga(1)$ is Zariski 
dense in $ SL (2,
\R ) \times  SL (2, \R )$ and $\A \supset k$. Or, one could also give
an elementary proof of the Zariski density in
$ GL (2,
\C ) \times  GL (2, \C )$ as follows.

Define  $\A^*_1 : = \overline{ \A^*} \cap ( GL (2, \C ) \times \{1\})$,
and $\A^*_2$ similarly (where $\overline{ \A^*}$ denotes the Zariski closure).

We observe that if $\A^*_1 = GL (2, \C )$, then  by extending to $\C$ 
the Galois
automorphism of $k$, we see that also $\A^*_2 = GL (2, \C )$,
thus there is nothing left to prove.

Else, both $\A^*_1 $ and $\A^*_2$ are proper algebraic subgroups
of  $ GL (2, \C )$, thus they are both solvable. Each  of the respective
projections  $ p_i (\overline{ \A^*} ) \subset  GL (2, \C )$
is surjective (by two reasons: if it would
yield a proper subgroup it follows that $\overline{ \A^*} $ is solvable,
a  contradiction; or, just use that $\A \otimes_j \C  = GL (2, \C)$).

Again   by extending to $\C$ the Galois
automorphism of $k$ we also see that $\overline{k}$
equals the centre $\C^* \times
\C^*$, and that $\overline{ \A^*}$ contains a commutative
subgroup of dimension equal to 4.

Thus, $ \overline{ \A^*} $ contains the direct product  $T_1 \times T_2$ of two
respective maximal tori of $ GL (2, \C )$. Since $ \overline{ \A^*} $ projects
onto $ GL (2, \C )$ by the first projection,
it contains also the union of the conjugates of   $T_1 \times \{ 1\}$,
thus, being closed, it contains $ GL (2, \C ) \times \{ 1\}$.
Similarly, it contains  $   \{ 1\} \times GL (2, \C )$ and we are done.

\QED

\begin{lem}
Let $\Delta \subset \PP SL (2, \R ) \times  \PP SL (2, \R )$
be a subgroup commensurable with $\Gamma$: then
"$\Delta \subset \A  $", more precisely $\Delta \subset (\PP j_1
\times\PP  j_2) (\A ) $.
\end{lem}

\Proof

Consider the inverse image $\De(1)$ of $\De$ inside
   $ SL (2, \R ) \times   SL (2, \R )$: then
$\De(1)$ is commensurable with $\Ga(1)$ and
there is a finite index subgroup $ \Gamma " \subset \Ga (1)$ such
that
each $\delta \in \Delta (1)$ normalizes
$ \Gamma ''  $.

Thus $\de$ normalizes the $k$-linear span of $\Ga ''$
inside $ M(2, \R) \times M (2, \R)$.
By the previous lemma a fortiori the $\Q$-linear span
of $\Ga ''$ equals
$\A$ .

   It follows then that $\de$ normalizes $\A$, and
by the Skolem-Noether Theorem
(cf.  e.g. \cite{blan}, Theorem III-4,  page 70)
it follows that there is an element $ \ga \in \A$ such that
conjugation of $\A$ by  $\de$ equals the inner automorphism
associated to $\ga$. Therefore we obtain that $\de \ga^{-1}$
centralizes $\A$.

Since however (cf. lemma above ) $\A$ is Zariski dense in
$ M (2, \R ) \times   M (2, \R )$, it follows
that the element $\de \ga^{-1}$ lies in the centre
   $ \{ \pm 1 \} \times \{ \pm 1\}$ of $ SL (2, \R ) \times   SL (2, \R )$,
whence the image of $\de$ inside
$\PP SL (2, \R ) \times  \PP SL (2, \R )$
lies in the image of $\A$.

\QED

\begin{lem}
Let $\de \in \A \setminus \{1\} $ yield a transformation
of $\HH^2$ which
has a fixed point : then the subfield $K_{\de} : = k [\de] \subset
\A$ is a cyclotomic extension $k [\zeta_m]$ where
$m \in \{3,4,5,6,8,10,12 \}$.
\end{lem}

The proof of the above lemma is contained in \cite{shav}, prop. 4.6,
and in the considerations following it. The main idea
is that if $\de$ has a fixed point, then it has finite order,
whence $K_{\de}$ is a cyclotomic extension: but
then the degree of the extension $ \Q \subset \Q [\zeta_m]$
divides 4, and one concludes calculating the $m$'s for which
the Euler function $\phi (m)$ divides 4.

\begin{df}\label{S}
Consider the greatest common multiple $120$ of the integers
appearing in the previous lemma, and let $ K$ be the
cyclotomic extension $k [\zeta_{120}]$.

For each intermediate field  $K'$ $ k \subset K' \subset K$
choose a prime ideal  $\PPP ' \subset \hol_k$ such that
   $\PPP '  \hol_{K'}$ is not primary. Such an ideal
exists for each such $K'$ by the density theorem (cf.  \cite{lang}.
VII, 4 , page 168)
and guarantees that $ k_{\PPP'} \otimes_k K'$ is
not a field (indeed , it is not an integral domain).

Let $\SSS ' : = \{  \PPP' \} \subset Spec (\hol_k)$,
and take $\SSS \subset Spec (\hol_k)$ as a  set of even
cardinality containing $ \SSS' $.

\end{df}

\begin{teo}
Let $k$ be a real quadratic field, and let $\A$ be the
indefinite division quaternion algebra corresponding,
by Hasse's theorem, to a choice of $\SSS$
made as in \ref{S}.

Then any subgroup  $\Delta \subset \PP SL (2, \R )
\times  \PP SL (2, \R )$
   commensurable with the subgroup $\Gamma$
associated to a maximal order $ \RR \subset \A$ acts freely
on $\HH^2$.
\end{teo}

\Proof
Assume that $\de \in \De$ is a nontrivial
element which does not act freely.
We have shown that $ \de \in \A$, and that
$K_{\de}$ is an intermediate field $K'$
between $k$ and $ K : = k [\zeta_{120}]$.

By our choice of $\PPP '$, it follows
that $\A \otimes_k k_{\PPP'}$ is a
division algebra; but on th other side
we have that $\A \otimes_k k_{\PPP'}$
contains $K_{\de} \otimes_k k_{\PPP'} =
K' \otimes_k k_{\PPP'} $ which is
not an integral domain. This is a contradiction.

\QED

Hence follows
\begin{teo}\label{quaternion}
Let $k$ be a real quadratic field, and let $\A$ be the
indefinite division quaternion algebra corresponding,
by Hasse's theorem, to a choice of $\SSS$
made as in \ref{S}.

Define $\FF_{\A}$ to be the family of subgroups
    $\Delta \subset \PP SL (2, \R ) \times  \PP SL (2, \R )$
   commensurable with a subgroup $\Gamma$
associated to a maximal order $ \RR \subset \A$.

Each $ \De \in \FF_{\A}$ acts freely
on $\HH^2$, and denote by
$S_{\De} : = \HH^2 / \De$
the corresponding algebraic surface.

Then the family of surfaces $ \{ S_{\De} | \De \in \FF_{\A} \}$
   is a union of Q.E.D. equivalence
classes.

\end{teo}

\Proof
Assume that a surface $S$ is Q.E.D. -equivalent to
$S_{\De}$: then it is Q.E.D. -equivalent to
$S_{\Ga}$, whence it corresponds to a subgroup
$\De ' $ commensurable with $\Ga$.

\QED

\begin{cor}\label{aleph}
There are infinitely many Q.E.D. -equivalence classes
of algebraic surfaces of general type.
\end{cor}

\Proof
It suffices to observe that the fundamental group
$\De $ of $S_{\De}$ has, by the cited theorem of Jost-Yau,
at most two embeddings inside
$  SL (2, \R ) \times    SL (2, \R )$
with isomorphic image acting freely and cocompactly.

These two are conjugate of each other, and for both
one sees by \ref{span} that the $\Q$-linear span of $\De$ equals
the  $\Q$-linear span of $\Ga$, which is indeed the
embedded quaternion algebra $\A$.

Since $\A$ determines its centre $k$ and its set of primes
$\SSS (\A)$, we see that to surfaces of the same Q.E.D. class
corresponds the same pair $ ( k, \SSS )$.

We saw however that the image is countable ( indeed,
for each choice of $k$, there are countably many choices
for $\SSS$), whence there are countably many Q.E.D. classes.

\QED

\section{Appendix by S\"onke Rollenske: Q.E.D. for Kodaira surfaces}

The aim of this appendix is to study the $QED$ equivalence
relation for Kodaira surfaces. More precisely we want to prove the following

\begin{teo}\label{kodaira}
Let $S$ be a minimal Kodaira surface. Then a smooth surface $S'$ is 
$QED$ equivalent to
$S$ if and only if $S'$ is itself a Kodaira surface. Thus Kodaira surfaces
constitute a single $QED$ equivalence class.
\end{teo}
The only if part of the teorem  is mostly an adaption of notes of
F. Catanese regarding the $QED$ equivalence for Hopf surfaces. Let us begin
with some preliminary considerations.

The surfaces of Kodaira dimension zero which are not K\"ahler are
called Kodaira surfaces. The minimal surfaces fall in the following two classes
of which Kodaira gave an explicit description (cf. \cite{kod1}, \cite{kod2}
and \cite{bpv}).

A minimal surface $S$ of Kodaira dimension zero is called a \emph{primary
Kodaira surface} if one of the following equivalent conditions holds:
\begin{itemize}
\item The first Betti number $b_1=3$.
\item $S$ is a holomorphic principal bundle of elliptic curves over an
elliptic curve, which is not topologically trivial.
\item $S$ is isomorphic to a Quotient $\C^2/G$ where $G$ is a group of affine
transformations generated by
\begin{gather*} g_i: (z_1,z_2)\mapsto (z_1+\alpha_i,z_2+\bar\alpha_i z_1
+\beta_i)\qquad i=1,\dots, 4
\intertext{with}\tag{$\ast$}\label{group}
  \alpha_1, \alpha_2=0\qquad\beta_1\bar\beta_2-\beta_2\bar\beta_1\neq 0\\
  \alpha_3\bar\alpha_4-\alpha_4\bar\alpha_3=m\beta_2\neq 0.
\end{gather*}
where $m$ is a positive integer.
The global holomorphic forms on $S$ are given by scalar multiples of (the
classes of) $dz_1$ and $dz_1\wedge dz_2$.
\end{itemize}

Sometimes a primary Kodaira surface admits a finite group of fixed point
free automorphisms such that the quotient is an elliptic quasi-bundle over
$\P^1$. Such a surface has $b_1=1$ and is called a \emph{secondary Kodaira
surface}.

A smooth surface $S$ of Kodaira dimension zero is a bimeromorphically 
equivalent
to a Kodaira surface if and only if $b_1(S)$ is  1 or 3. From this we get
immediately the follwing

\begin{cor}\label{defKodaira}
Assume that we have a flat family $F : \X \ra \Delta$ over the unit
disk with special fibre
a compact complex surface $X$ with canonical singularities, and with
another fibre which is a smooth Kodaira surface $S$.
Then $X$ is bimeromorphic to a Kodaira surface.
\end{cor}
\Proof
By Tyurina's result (\cite{tyu}) the minimal resolution $Z$ of the
singularities of $X$
is a surface diffeomorphic to $S$ and so has uneven Betti number. By
Theorem S7 of \cite{f-m} (p. 224) it has also Kodaira dimension zero and thus
$Z$ is a
Kodaira surface. \QED

Now let us analyse the automorphisms of minimal Kodaira surfaces by making use
of the
above description. Here we follow closely Kodaira:
Let $\Gamma$ be a finite group of automorphisms of $S=\C/G$. By pulling back to
the universal covering we get an extension of finite index
\[1\to G\to \Gamma'\to\Gamma\to 1\]
where $\Gamma'$ is a group of automorphisms of $\C^2$. Now let $\phi=(\phi_1,
\phi_2)$ be in $\Gamma'$. The linear action of $\Gamma$ on $H^0(S,\Omega_S^p)$
becomes
\begin{gather*} \phi^*dz_1=d\phi_1=\sigma dz_1 \\
\phi^*(dz_1\wedge dz_2)=d\phi_1\wedge d\phi_2=\kappa dz_1\wedge dz_2 =\sigma
dz_1\wedge d\phi_2
\end{gather*}
with $\sigma, \kappa\in \C$ and consequently there exist a function $h(z_1)$
and a constant $h_0$ such that \[ \phi_1= \sigma z_1+ h_0 \text{ and } \phi_2
=\frac{\kappa}{\sigma}z_2 + h(z_1). \]
Since $G$ is a normal subgroup of finite index in $\Gamma'$ we have $\phi^n\in
G$ so ${\sigma}$ is a root of unity and for every generator of $G$
there exists an element $\tilde
g_i(z_1,z_2)=(z_1+a_i,z_2+\bar{a}_i z_1
+ b_i) \in G$ such that $\phi \circ g_i= \tilde g_i\circ \phi$.
Calculating both sides we get $\sigma \alpha_i=a_i$ and also
- after deriving the second part with respect to $z_1$:
\[
h'(z_1)-h'(z_1+\alpha_i)=\bar\alpha_i\left(\frac{\kappa}{\sigma}
-\sigma\bar\sigma\right)=\bar\alpha_i\left(\frac{\kappa}{\sigma}
-1\right).\]
Now $h''(z_1)$ is constant because it has two linear independent periods
$\alpha_3, \alpha_4$ and consequently
\begin{gather*}h''(z_1) \alpha_i=\bar\alpha_i\left(\frac{\kappa}{\sigma}
-1\right)\qquad i=3,4.\\ \Rightarrow h''=\left(\frac{\kappa}{\sigma}
-1\right)=0.
\end{gather*}
Hence we have
\begin{equation}\label{autom}
\phi(z_1, z_2)=\begin{pmatrix}\sigma & h_1\\0 &1\end{pmatrix}
\binom{z_1}{z_2}+
\binom{h_0}{h_2}\end{equation}

Now assume that an automorphism $\bar\phi\in \Gamma$ has fixed points.
We take a lift $\phi\in\Gamma'$ and can assume (by multiplying by an element of
$G$ if
necessary) that $\phi$ itself has a fixed
point. By
(\ref{autom}) this is the case iff the equation
\[0= \left[ \begin{pmatrix}\sigma & h_1\\0 & 1\end{pmatrix}-
\begin{pmatrix}1 &0\\0&1\end{pmatrix}\right]
\binom{z_1}{z_2}+
\binom{\gamma}{h_2}=\begin{pmatrix}\sigma-1 & h_1\\0 &0\end{pmatrix}
\binom{z_1}{z_2}+
\binom{\gamma}{h_2}\]
has a solution. Clearly the same arguments work for secondary Kodaira surfaces.

Since an automorphism maps rational curves to rational curves and therefore
covers a unique automorphism of the corresponding minimal model we have shown
the first part of the following
\begin{prop}\label{autoKodaira} If an automorphism of finite order of a
Kodaira surface has fixed points, it has fixed points in codimension one. In
particular if $S$
is a Kodaira surface and $f:S\to X$ is a quasi-\'etale map where $X$ has
canonical singularities, than than $f$ is \'etale and $X$ is in fact a
smooth Kodaira surface.\end{prop}
\Proof
By corollary \ref{groupquotient} there is a Kodaira surface $W$ and a finite
group $G$ acting freely in codimension 1 on $W$ such that $X$ is a
birational image of $W/G$ by a small contraction. But the first part implies
that the action of $G$ is free, hence the quotient is smooth, $W/G\cong X$ and
$f$ is itself \'etale.\QED

\begin{lem}\label{upKodaira}
Assume that  $f : X \to Y$ is a quasi \'etale morphism where $Y$ has
canonical singularities
and is bimeromorphic to a smooth Kodaira surface $S$.
Then $X$ is bimeromorphic to a Kodaira surface.
\end{lem}
\Proof
Without loss of generality, we may assume that $ \pi : S \ra Y$ be the minimal
resolution of the singularities of $Y$, and that $S$ be the blow up
$ p : S \ra Z$ of a minimal Kodaira surface $Z$.

By definition, there are finite sets $\Sigma_Y \subset Y, \Sigma_X \subset X $
such that $ X - \Sigma_X \ra Y - \Sigma_Y$ is a finite unramified covering.

By pull back, we obtain a finite unramified covering of $ S -
\pi^{-1} (\Sigma_Y)$.
Now, $\pi^{-1} (\Sigma_Y)$ consists of a finite set plus a finite union of
smooth rational curves with self intersection $-2$.

But there is no rational curve on a  minimal Kodaira surface $Z$, the
universal covering
    being $\C^2 $,  whence $\pi^{-1} (\Sigma_Y)$ maps onto a
finite set on $Z$,
and $X$ is bimeromorphic to a quasi \'etale covering $W$ of $Z$.
Since $Z$ is smooth, $W$ is a finite unramified cover of $Z$; in particular it
is not K\"ahler.

Every minimal Kodaira surface admits a volume-preserving complex structure and
so this is also true for $W$.  Kodaira classified these surfaces completely in
(\cite{kod2}, Theorem 39) and by  his results $W$ is covered either by a
K3 surface, a complex torus or a primary Kodaira surface. Since $W$ covers $S$
only the third case occurs and $W$ is itself a Kodaira surface.\QED

{\itshape Proof of the Theorem.}
First let $S$ be a minimal Kodaira surface. We want to show, that 
every surface $QED$ equivalent to $S$ is a Kodaira surface.
By \ref{defKodaira} and \ref{upKodaira} it suffices to show that if we have
a quasi - \'etale morphism $p : X \ra Y$ where $X, Y$ have canonical
singularities
and $X$ is bimeromorphic to $S$, then $Y$ is also bimeromorphic to a
Kodaira surface.

By taking the normal closure, and applying \ref{upKodaira} we may assume that
$p : X \ra Y$ is the quotient map by the action of a finite group $G$.
Now  $Y$ is bimeromorphic to $S/G$ which is a Kodaira surface by 
Proposition \ref{autoKodaira}.

It remains to show, that all Kodaira surfaces are $QED$ equivalent. Let
$S_0=\C^2/G_0$ be the primary Kodaira surface given by (\ref{group}) with
$\beta_1=\alpha_3=1$, $\beta_2=2\alpha_4=2i$ and $\beta_3=\beta_4=0$.
We have the relation $\alpha_3\bar\alpha_4-\alpha_4\bar\alpha_3=\beta_2$ and
thus the fundamental group of $S_0$ is isomorphic as an abstract group to $F/R$
where $F$ is the free group on generators $f_1,\dots, f_4$ and $R$ is the
subgroup generated by the relations
\[ [f_i,f_j]=\begin{cases} f_2 & i=3,j=4\\ 0 & 
\text{otherwise}\end{cases} \qquad 1\leq i<j\leq 4.\]

It clearly suffices to show that every primary Kodaira surface is 
$QED$ equivalent to $S_0$. First consider an arbitrary $S=\C^2/G$ 
with $G$ as in (\ref{group}). By changing $\alpha_3$ to 
$\alpha_3'=\frac{\alpha_3}{m}$ we get another group $G'$, a surface 
$S'=\C^2/{G'}$ and finite covering maps
\[\xymatrix{ S \ar[r]\ar[d] & \C/{<\alpha_3,\alpha_4>}\ar[d]\\
  S'\ar[r] &\C/{<\alpha_3',\alpha_4>}}\]
We have the relation
$\alpha_3'\bar\alpha_4-\alpha_4\bar\alpha_3'=\beta_2$ and thus 
$\pi_1(S')=G'$ is
isomorphic to $F/R$. By Corollary II.7.17 of \cite{f-m} it follows that $S$ is
deformation equivalent to $S_0$ or to $S_0^{conj}=\C^2/{\bar G_0}$ with the
conjugated complex structure. But an easy calulation shows that $G_0=\bar G_0$
thus $S_0=S_0^{conj}$ and consequently $S$ is $QED$ equivalent to $S_0$ which
concludes the proof. \QED \\
The last part of the proof can also be obtained using the description of the
moduli space obtained by Borcea in \cite{borc}.

\section{Final remarks}

\noindent

{\bf Note 1.}
As remarked by Frederic Campana in the footnote to \cite{cam},
the equivalence
relation introduced by him is only apparently similar to ours, but
indeed quite different,
cf. section 5 of \cite{cam} and our main theorems for special surfaces.

{\bf Note 2.}
Claire Voisin pointed out  that  the  decision to use also the notion of
K-equivalence (introduced in \cite{voi3}) might  lead to other interesting
equivalence relations preserving the Kodaira dimension.

{\bf Acknowledgement} The main results on $\C$-Q.E.D. of K\"ahler surfaces
were announced in the Luminy G.A.C. Conference in december 2001, and at other
conferences (cf. the Abstracts of the Fano Conference, Torino october 2002).
I thank G. Laumon for convincing me of the importance to check also
A.Q.E.D. equivalence, Thomas Peternell for a  useful suggestion
concerning proposition \ref{peano},
Igor Dolgachev for a useful discussion on logarithmic transformations,
Ciro Ciliberto for a useful discussion, 
and especially Fritz Grunewald and S\"onke Rollenske for their contribution and for
several critical remarks.

The research of the  author was performed in the realm  of the
    SCHWERPUNKT "Globale Methoden in der komplexen Geometrie",
and of the EAGER EEC Project.

\vfill

\bigskip

\noindent
Prof. Fabrizio Catanese\\
Lehrstuhl Mathematik VIII\\
Universit\"at Bayreuth, NWII\\
    D-95440 Bayreuth, Germany

e-mail: Fabrizio.Catanese@uni-bayreuth.de

\end{document}